\documentclass{amsart} 

\usepackage{amsmath, amssymb, latexsym, amscd, mathabx}
\usepackage{exscale, cite, eps fig, graphics}
\usepackage{comment}
\usepackage[colorlinks]{hyperref}

\renewcommand{\geq}{\geqslant}

\renewcommand{\l}{\langle}
\renewcommand{\r}{\rangle}

\renewcommand{\c}{c_{\textrm{\tiny ww}}}
\newcommand{\cc}{c_{\emph{\tiny ww}}}
\renewcommand{\k}{\kappa}

\renewcommand{\L}{\mathcal{L}}

\newcommand{\p}{\boldsymbol{\phi}}

\renewcommand{\v}{\mathbf{v}}

\newtheorem*{main-theorem}{Main Theorem}
\newtheorem*{remark*}{Remark}
\numberwithin{equation}{section}

\usepackage{subcaption}

\usepackage{scalerel,stackengine}
\stackMath
\newcommand\reallywidehat[1]{%
\savestack{\tmpbox}{\stretchto{%
  \scaleto{%
    \scalerel*[\widthof{\ensuremath{#1}}]{\kern-.6pt\bigwedge\kern-.6pt}%
    {\rule[-\textheight/2]{1ex}{\textheight}}%WIDTH-LIMITED BIG WEDGE
  }{\textheight}% 
}{0.5ex}}%
\stackon[1pt]{#1}{\tmpbox}%
}
\title[Modulational instability in full-dispersion shallow water models]{Comparison of modulational instabilities in full-dispersion shallow water models}

\author[Pandey]{Ashish~Kumar~Pandey}
\address{Department of Mathematics, University of Illinois at Urbana-Champaign, Urbana, IL 61801 USA}
\email{akpande2@illinois.edu}  

\date{\today}

\begin{document}

\maketitle

\begin{abstract}
We study the modulational instability of a shallow water model, with and without surface tension, which generalizes the Whitham equation to include bi-directional propagation. Without surface tension, the small amplitude periodic traveling waves are modulationally unstable if their wave number is greater than a critical wave number predicting a Benjamin-Feir type instability and the result qualitatively agrees with the shallow water model in \cite{HP2}. With surface tension, the result qualitatively agrees with the physical problem except for the large surface tension limit which is accurately predicted by the shallow water model in \cite{HP2}. We also compare the results with the Whitham and full-dispersion Camassa-Holm equations.
\end{abstract}

%%%%%%%%%%%%%%%%%%%%%%%%%%%%%%%%%%%%%%%%%%%%%%%%%%
%%%%%%%%%%%%%%%%%%%%%%%%%%%%%%%%%%%%%%%%%%%%%%%%%%
%%%%%%%%%%%%%%%%%%%%%%%%%%%%%%%%%%%%%%%%%%%%%%%%%%
\section{Introduction}\label{sec:intro}
%%%%%%%%%%%%%%%%%%%%%%%%%%%%%%%%%%%%%%%%%%%%%%%%%%
%%%%%%%%%%%%%%%%%%%%%%%%%%%%%%%%%%%%%%%%%%%%%%%%%%
%%%%%%%%%%%%%%%%%%%%%%%%%%%%%%%%%%%%%%%%%%%%%%%%%%

In a laboratory experiment, Benjamin and Feir \cite{BF} discovered that a slight modulation imposed on nonlinear Stokes waves on the water surface results into the disintegration of waves over time and also explained it theoretically for deep water. Whitham \cite{Whitham1967} using formal arguments showed that the Stokes waves are unstable when $\k h > 1.363\dots$, where $\k$ is the wave number and $h$ is the undisturbed fluid depth. This instability is called {\em modulational instability} or {\em Benjamin-Feir instability} in the context of water waves. The phenomenon was also described by formally deriving a nonlinear Schr\"{o}dinger equation for the envelope of the Stokes waves, see \cite{BN, Zakharov-WW,HO}. A rigorous proof of the Benjamin-Feir instability for the finite depth case has been given in \cite{BM1995}. The modulational instability in capillary-gravity waves has also been investigated, see\cite{Kawahara} and \cite{DR}, for instance.           

From the condition $\k h > 1.363\dots$, it can be inferred that the Benjamin-Feir instability is a high-frequency phenomenon and may not be seen in water wave models which approximate full water wave equations in the long-wave regime. This manifests in the Korteweg-de Vries (KdV) equation 
\begin{equation*}\label{E:KdV}
u_t +u_x + \frac16 u_{xxx} + (u^2)_x = 0,
\end{equation*}
for which all periodic traveling waves are modulationally stable; see \cite{BHJ} and references therein. Here, $u(x,t)$ is a real-valued function describing the water velocity at position $x$ and time $t$ from the undisturbed fluid depth. Note that the dispersion (which we describe by the phase velocity of the linear part) $c_{\textrm{\tiny KdV}}(\k)$ of the KdV equation agrees with the dispersion of full water wave equations 
\begin{equation}\label{def:c}
\c(\k) := \sqrt{\frac{\tanh(\k)}{\k}} = 1-\frac16 \k^2 + O(\k^4) =c_{\textrm{\tiny KdV}}(\k) + O(\k^4),
\end{equation}
only for small frequencies. In the presence of surface tension (see \cite{Whitham}, for instance),
\begin{equation}\label{def:c1}
\c(\k)=\sqrt{(1+T \k^2)\frac{\tanh(\k)}{\k}}
\end{equation}
replaces \eqref{def:c}, where $T\geq 0$ is the coefficient of surface tension. For $T=0$, \eqref{def:c1} reduces to \eqref{def:c}. In fact, the full-dispersion of water waves has two branches, a positive branch, $\c(\k)$, and a negative branch, $-\c(\k)$. A water wave model with dispersion same as the full water wave equations is called a {\em full dispersion water wave model}.

The Whitham equation   
\begin{equation}\label{E:Whitham}
\partial_t u+\c(|\partial_x|)\partial_x u
+u\partial_x u=0, 
\end{equation}
generalizes the KdV equation to include full-dispersion (or more precisely, the positive branch, $\c(\k)$) of water wave problem. Here, $\c(|\partial_x|)$ is Fourier multiplier with symbol $\c(\k)$, 
\begin{equation*}
\reallywidehat{\c(|\partial_x|) v}(\k)= \c(\k) \widehat{v}(\k).
\end{equation*}
The Whitham equation was introduced by Whitham \cite{Whitham} as a model to study breaking of nonlinear dispersive water waves. The Whitham's conjecture regarding wave breaking in \eqref{E:Whitham} and \eqref{def:c} was settled in \cite{Hur-breaking1}. In \cite{HJ2}, small amplitude periodic traveling waves of \eqref{E:Whitham} and \eqref{def:c} have been shown to display Benjamin-Feir instability with a critical wave number, $1.146\dots$. The effects of surface tension on the modulational instability in \eqref{E:Whitham} and \eqref{def:c1} were examined in \cite{HJ3} and it explained the effects in accordance with the physical problem \cite{DR} except for the large surface tension limit. For the full water wave problem, for large surface tensions, the stability changes to instability about a critical wave number similar to the Benjamin-Feir instability and this critical wave number converges in the large surface tension limit \cite{DR}. The Whitham equation explains this phenomenon for large surface tensions, but the limit of critical wave number diverges in the large surface tension limit \cite{HJ3}. Undoubtedly, the Whitham equation improves upon the KdV equation in explaining phenomena of water waves like wave breaking and Benjamin-Feir instability but it has its shortcomings, for example, it fails to explain the modulational instability accurately in the large surface tension limit.

A full-dispersion shallow water (FDSW-I) model generalizing the nonlinear shallow water equations
\begin{equation*}\label{E:shallow}
\begin{aligned}
&\partial_t\eta+\partial_x(u(1+\eta))=0,\\
&\partial_tu+\partial_x\eta+u\partial_xu=0
\end{aligned}
\end{equation*}
to include full-dispersion of water waves has been proposed in \cite{HT2,HP2} given by
\begin{equation}\label{E:BW}
\begin{aligned}
&\partial_t\eta+\partial_x(u(1+\eta))= 0,\\
&\partial_tu+\c^2(|\partial_x|)\partial_x\eta+ u\partial_xu=0,
\end{aligned}
\end{equation}
where $\c^2(|\partial_x|)$ is Fourier multiplier with symbol $\c^2(\k)$, 
\begin{equation*}
\reallywidehat{\c^2(|\partial_x|) v}(\k)= \c^2(\k) \widehat{v}(\k).
\end{equation*}
Here, $\eta(x,t)$ represents the surface displacement. The wave breaking of \eqref{E:BW} and \eqref{def:c} has also been settled in \cite{HT2}. The small amplitude periodic traveling waves of \eqref{E:BW} and \eqref{def:c} exhibit the Benjamin-Feir instability with a critical wave number, $1.610\dots$ \cite{HP2}. Moreover, the effects of surface tension on modulational instability in \eqref{E:BW} and \eqref{def:c1} turned out to be qualitatively same with the physical problem \cite{Kawahara} for all values of surface tension \cite{HP2}. Therefore, the FDSW-I model explains water waves' phenomena like wave breaking and Benjamin-Feir instability similar to the Whitham equation. In addition to being a bi-directional model, the FDSW-I model overcomes the shortcoming of the Whitham equation in explaining the effects of surface tension on modulational instability in the large surface tension limit.   

Another full-dispersion shallow water (FDSW-II) model can be obtained by shifting $\c^2(|\partial_x|)$ in the first equation of \eqref{E:BW}
\begin{equation}\label{E:main}
\begin{aligned}
&\partial_t\eta+\partial_x(\c^2(|\partial_x|) u+u\eta)=0,\\
&\partial_tu+\partial_x\eta+u\partial_x u=0.
\end{aligned}
\end{equation}
In fact, this model was proposed in \cite{MKD}. In Section~\ref{sec:MI}, we take matters further for \eqref{E:main} and \eqref{def:c1} and derive an index \eqref{def:ind} which if negative provides modulational instability. The index \eqref{def:ind} is explicitly given as a function of the wave number $\k$ and the coefficient of surface tension $T$. A Benjamin-Feir type instability with a critical wave number $1.008\dots$ has been obtained by examining the sign of the index for $T=0$. To study the effects of surface tension on modulational instability, we inspect the sign of the index for $T>0$. In the $\k$ and $\k \sqrt{T}$ plane, we determine the regions of stability and instability; see Figure~\ref{fig1:d} for details. The result qualitatively agrees with the physical problem \cite{Kawahara} except for the large surface tension limit.

A unidirectional model which overcomes the failure of the Whitham equation in explaining the effects of surface tension on modulational instability in the large surface tension limit is the full-dispersion Camassa-Holm (FDCH) equation
\begin{equation}\label{E:FDCH}
\eta_t+c_{\rm ww}(|\partial_x|)\eta_x+\frac{3\eta}{1+\sqrt{1+\eta}}\eta_x
=-\Big(\frac{5}{12}\eta\eta_{xxx}+\frac{23}{24}\eta_x\eta_{xx}\Big).
\end{equation}
The FDCH equation has been proposed in \cite{Lannes} and its modulational instability, with or without surface tension, has been examined in \cite{HP3}. Small periodic traveling waves of \eqref{E:FDCH} and \eqref{def:c} show Benjamin-Feir instability with a critical wave number, $1.42\dots$ \cite{HP3}. Moreover, the effects of surface tension on the modulational instability in \eqref{E:FDCH} and \eqref{def:c1} qualitatively agree with the physical problem \cite{HP3}.

In Section~\ref{sec:T>0}, we compare results on modulational instability in the Whitham, FDSW-I, FDSW-II, and FDCH equations. In the absence of surface tension, we find that all these models exhibit a Benjamin-Feir type instability. With surface tension, the results qualitatively agree with the physical problem except for the large surface tension limit. The FDSW-I and FDCH equations accurately predict the large surface tension limit of the critical wave number while the Whitham and FDSW-II equations fail to do so.              

%%%%%%%%%%%%%%%%%%%%%%%%%%%%%%%%%%%%%%%%%%%%%%%%%%
%%%%%%%%%%%%%%%%%%%%%%%%%%%%%%%%%%%%%%%%%%%%%%%%%%
%%%%%%%%%%%%%%%%%%%%%%%%%%%%%%%%%%%%%%%%%%%%%%%%%%
\section{Modulational instability in \eqref{E:main} and \eqref{def:c1}}\label{sec:MI}
%%%%%%%%%%%%%%%%%%%%%%%%%%%%%%%%%%%%%%%%%%%%%%%%%%
%%%%%%%%%%%%%%%%%%%%%%%%%%%%%%%%%%%%%%%%%%%%%%%%%%
%%%%%%%%%%%%%%%%%%%%%%%%%%%%%%%%%%%%%%%%%%%%%%%%%%

We derive a modulational instability index for \eqref{E:main} and \eqref{def:c1}. The details of the proof follow along Section~$2$ and $3$ of \cite{HP2}. We briefly present the main ideas and record relevant expressions. 

By a traveling wave of \eqref{E:main}, we mean a stationary solution of form $(\eta,u)(x,t) = (\eta,u)(x-ct)$ for some $c>0$. Further, we take $\eta$ and $u$ to be $2\pi$-periodic functions of $z=\k x$. The result becomes, by quadrature, 
\begin{equation}\label{E:periodic}
\begin{aligned}
&-c\eta+\c^2(\k|\partial_z|)u+u\eta=0,\\
&-cu+\eta+\frac12u^2=0.
\end{aligned}
\end{equation}

For every $T\geq 0$, the existence of a smooth $\eta$ and $u$ satisfying \eqref{E:periodic} follows from a Lyapunov-Schmidt procedure with the condition that the kernel of the linearized operator of \eqref{E:periodic} is only one dimensional in the space of even functions in $H^1(\mathbb{T})$(see, \cite[Section~2]{HP2}, for details). The aforementioned kernel is one dimensional if $\c(\k)$ is monotone in $\k$ which is the case for all $T\geq 0$ except for $0<T<1/3$ (see \cite[Figure~7]{HP2}). For $0<T<1/3$, the kernel will be one dimensional if
\begin{equation}\label{E:k-condition}
\c(\k)\neq \c(n\k), \qquad n=2,3,\dots.
\end{equation}
For any $T\geq 0$ and $\k>0$ satisfying \eqref{E:k-condition}, existence of a one parameter family of smooth solutions of \eqref{E:periodic} can be proved following the proof of Theorem~$2.1$ in \cite{HP2}. The small amplitude expansion of solutions is given as 
\begin{subequations}\label{E:hu-small}
\begin{align}
\eta(a;\k)(z)=&a\c(\k) \cos z \notag+a^2\Big(\c(\k)h_0-\frac14+\Big(\c(\k)h_2-\frac14\Big)\cos2z\Big )+O(a^3), \notag \\
u(a;\k)(z)=&a\cos z +a^2(h_0+h_2\cos 2z)+O(a^3),\notag
\intertext{and}
c(a;\k)=&\c(\k)+\frac32 a^2 \Big(h_0+\frac12 h_2-\frac{1}{8\c(\k)}\Big)+O(a^3)\hspace*{-25pt}\notag
\end{align}
\end{subequations}
as $a \to 0$, where
\begin{equation}\label{def:h02}
h_0=\frac34\frac{\c(\k)}{\c^2(\k)-1}\quad\text{and}\quad 
h_2=\frac34\frac{\c(\k)}{\c^2(\k)-\c^2(2\k)}.
\end{equation}

We linearize \eqref{E:main} about $\eta$ and $u$ in the coordinate frame moving at the speed $c>0$ and seek a solution of the form $\v(z,t)=e^{\lambda \k t}\v(z)$, $\lambda\in\mathbb{C}$, to arrive at
\begin{equation*}\label{def:L}
\lambda\v=\partial_z\begin{pmatrix} c-u &-\c^2(\k|\partial_z|)-\eta \\ -1 & c-u \end{pmatrix}\v
=:\mathcal{L}(a;\k)\v.
\end{equation*}
Using Floquet theory (see \cite{BHJ,HP2}, for details), the  $L^2(\mathbb{R}) \times L^2(\mathbb{R})$ spectrum of $\L$ is decomposed into $L^2(\mathbb{T})\times L^2(\mathbb{T})$ spectra of $\L(\xi,a)$'s for $\xi\in(-1/2,1/2]$ defined by
\[
\L(\xi,a)\v(\xi):=e^{-i\xi z}\L e^{i\xi z}\v(\xi).
\]
For any $\xi \in(-1/2,1/2]$,  the $L^2(\mathbb{T})\times L^2(\mathbb{T})$ spectrum of $\L(\xi,a)$ consists of eigenvalues with finite multiplicities. A straightforward calculation shows that zero is an eigenvalue of $\L(0,0)$ with multiplicity four. For $|a|\neq 0$, zero continues to be an eigenvalue of $\L(0,a)$ with a four-dimensional generalized eigenspace. For $|\xi|$ and $|a|$ small, we are interested in the eigenvalues of $\L(\xi,a)$ bifurcating from the zero eigenvalue of $\L(0,a)$. For this purpose, we extend the eigenspace for the zero eigenvalue of $\L(0,a)$ to construct a four-dimensional eigenspace for the bifurcating eigenvalues of $\L(\xi,a)$, for $|\xi|$ and $|a|$ small. This eigenspace is spanned by (see \cite[Lemma~3.1]{HP2})  
\begin{equation}\label{def:p}
\begin{aligned}
\p_1(\xi,a)(z)=&\begin{pmatrix}\c(\k) \\ 1 \end{pmatrix}\cos z 
+i\xi \frac{\k\c'(\k)}{\c^2(\k)+1}\begin{pmatrix} 1\\ -\c(\k)\end{pmatrix}\sin z \\
&+\frac{a}{4\c(\k)}\begin{pmatrix} -\c(\k)(1+4\c(\k)h_2) \\ 1-4\c(\k)h_2\end{pmatrix} \\
&+\frac{a}{2}\begin{pmatrix}4\c(k)h_2-1\\ 4h_2\end{pmatrix}\cos 2z+\xi^2 \mathbf{p}_2\cos z+O(\xi^3+\xi^2a+a^2), \\
\p_2(\xi,a)(z)=&\begin{pmatrix}\c(\k) \\ 1 \end{pmatrix}\sin z 
-i\xi \frac{\k\c'(\k)}{\c^2(\k)+1}\begin{pmatrix} 1\\ -\c(\k)\end{pmatrix}\cos z\\
&+\frac{a}{2}\begin{pmatrix}4\c(k)h_2-1\\ 4h_2\end{pmatrix}\sin 2z
+\xi^2 \mathbf{p}_2\sin z+O(\xi^3+\xi^2a+a^2),\\
\p_3(\xi,a)(z)=&\begin{pmatrix} 2\c(\k) \\ -1\end{pmatrix}+a\begin{pmatrix}1\\ 0\end{pmatrix}\cos z
-\frac16\xi^2\k^2\c(\k)\begin{pmatrix} 1\\0\end{pmatrix}+O(\xi^3+\xi^2a+a^2), \\
\p_4(\xi,a)(z)=&\begin{pmatrix} 1\\2\c(\k)\end{pmatrix}+\frac{a}{2\c(\k)}\begin{pmatrix}1\\ 0 \end{pmatrix}\cos z
-\frac{1}{12}\xi^2 \k^2\begin{pmatrix} 1\\0\end{pmatrix}+O(\xi^3+\xi^2a+a^2)
\end{aligned}
\end{equation}
up to orders of $\xi^2$ and $a$ as $\xi$, $a\to 0$, where $h_2$ is in \eqref{def:h02} and 
\begin{equation*}\label{def:p2}
\mathbf{p}_2=\frac12\frac{\k^2}{\c^2(\k)+1} 
\begin{pmatrix}{\displaystyle -3\frac{(\c(\c')^2)(\k)}{\c^2(\k)+1}+\c''(\k)} \\ {\displaystyle \c'(\k)^2\frac{2\c^2(\k)-1}{\c^2(\k)+1}-(\c\c'')(\k)}
 \end{pmatrix}.
\end{equation*}

For $|\xi|$ and $|a|$ small, the four eigenvalues of $\L(\xi,a)$ bifurcating from zero eigenvalue coincide with the roots of $\det(\mathbf{L}-\lambda\mathbf{I})$ up to orders of $\xi^2$ and $a$ (see \cite[Section~4.3.5]{K}, for instance, for details), where
\begin{equation}\label{def:L3}
\mathbf{L}(\xi,a)=\left(\frac{\l\L(\xi,a)\p_k(\xi,a),\p_\ell(\xi,a)\r}{\l\p_k(\xi,a),\p_k(\xi,a)\r}\right)_{k,\ell=1,2,3,4}
\end{equation}
and 
\begin{equation}\label{def:I3}
\mathbf{I}(\xi,a)=\left(\frac{\l\p_k(\xi,a),\p_\ell(\xi,a)\r}{\l\p_k(\xi,a),\p_k(\xi,a)\r}\right)_{k,\ell=1,2,3,4},
\end{equation}
where $\p_1$, $\p_2$, $\p_3$, $\p_4$ are in \eqref{def:p} and $\langle\,,\rangle$ means the $L^2(\mathbb{T})\times L^2(\mathbb{T})$ inner product. This amounts to the fact that restricted on a four-dimensional eigenspace, $\L(\xi,a)$ can be defined by the $4\times 4$ matrix $\mathbf{L}(\xi,a)$ obtained by calculating its action on the basis $\{\p_1, \p_2, \p_3, \p_4\}$. Therefore, the eigenvalues of the resulting matrix are given by the roots of its characteristic polynomial $\det(\mathbf{L}-\lambda\mathbf{I})$, where $\mathbf{I}$ is the projection of the identity onto the eigenspace.

We omit all the details of the calculation as it is very similar to \cite{HP2} and report that \eqref{def:L3} becomes
\begin{equation}\label{E:L3}
\begin{split}
\mathbf{L}(\xi,a)=
&\frac14 a(\c^2(\k)+1)\begin{pmatrix} 0&0&0&0\\0&0&0&0\\0&0&0&0\\0&1&0&0\end{pmatrix}\\
&+i\xi\begin{pmatrix} -\k\c'(\k)&0&0&0\\ 0&-\k\c'(\k)&0&0\\ 
0&0&\c(\k)\frac{4\c^2(\k)+5}{4\c^2(\k)+1}&-\frac{4\c^2(\k)-1}{4\c^2(\k)+1}\\
0&0&-\frac{4\c^2(\k)-1}{4\c^2(\k)+1}&\c(\k)\frac{4\c^2(\k)-3}{4\c^2(\k)+1}\end{pmatrix} \\
&+i\xi a\,L\begin{pmatrix} 0&0&2\c(\k)&1\\ 0&0&0&0\\0&0&0&0\\0&0&0&0 \end{pmatrix} 
+i\xi a\frac{1}{2(4\c^2(\k)+1)} \begin{pmatrix}0&0&0&0\\0&0&0&0\\
L_{31}&0&0&0\\
L_{41}&0&0&0\end{pmatrix}\\
&+\frac12\xi^2\k(2\c'(\k)+\k\c''(\k))\begin{pmatrix}0&1&0&0\\ -1&0&0&0\\ 0&0&0&0\\0&0&0&0 \end{pmatrix}
+O(\xi^3+\xi^2 a+a^2)
\end{split}
\end{equation}
as $\xi$, $a\to0$, where 
\begin{align*}
L=&\frac{1}{2\c(\k)(\c^2(\k)+1)}(4\c(\k)(1-\c^2(\k))h_2-\k \c(\k)\c'(\k)-1-5\c^2(\k)), \\
L_{31}=&4\c(\k)(1-\c^2(\k))h_2-2-\c^2(\k),\\
L_{41}=&\frac{1}{2\c(\k)}(4\c(\k)(1-\c^2(\k))h_2-2(\c^2(\k)+1)-4\c^4(\k)),
\end{align*}
and $h_2$ is in \eqref{def:h02}. Moreover, \eqref{def:I3} becomes
\begin{equation}\label{E:I3}
\begin{split}
\mathbf{I}(\xi,a) =&\mathbf{I}+a\frac{4\c(\k)(1-2\c^2(\k))h_2-1}{4\c(\k)(\c^2(\k)+1)(4\c^2(\k)+1)}
\begin{pmatrix}0&0& 2(4\c^2(\k)+1) &0\\
0&0&0&0\\\c^2(\k)+1&0&0&0\\0&0&0&0 \end{pmatrix} \\
&+a\frac{1-6\c(\k)h_2}{2(\c^2(\k)+1)(4\c^2(\k)+1)}\begin{pmatrix}0&0&0&2(4\c^2(\k)+1)\\0&0&0&0\\
0&0&0&0\\
\c^2(\k)+1&0&0&0 
\end{pmatrix} \\
&-i\xi a\frac{\k \c'(\k)}{4\c(\k)(\c^2(\k)+1)^2(4\c^2(\k)+1)}\\
&\qquad \qquad \quad\begin{pmatrix}0&0&0&0\\
0&0&4\c(\k)(4\c^2(\k)+1)&2(4\c^2(\k)+1)\\
0&2\c(\k)(\c^2(\k)+1)&0&0\\
0&\c^2(\k)+1&0&0\end{pmatrix}\\
&+O(\xi^3+\xi^2 a+a^2) 
\end{split}
\end{equation}
as $\xi$, $a\to 0$, where $\mathbf{I}$ is the $4\times4$ identity matrix. The analysis of the roots of quartic polynomial, $\det(\mathbf{L}-\lambda\mathbf{I})$, in $\lambda$ can be analyzed using discriminants (see \cite[Section~3.7]{HP2}, for details), and we derive a modulational instability index for \eqref{E:main} and \eqref{def:c1} given by 
\begin{equation} \label{def:ind}
\Delta(\k):=\frac{i_1(\k)i_2(\k)}{i_3(\k)}i_4(\k),
\end{equation}
where 
\begin{align}
i_1(\k)=&(\k \c(\k))'', \notag \\
i_2(\k)=&((\k \c(\k))')^2-1, \notag \\
i_3(\k)=&\c^2(\k)-\c^2(2\k) \notag,
\intertext{and}
i_4(\k)=&9\c^2(\k)i_2(\k) \label{def:i4}+i_3(\k)(3+15\c^2(\k)+6\k\c(\k)\cc'(\k)-\k^2(\c'(\k))^2).
\end{align}
A sufficiently small, $2\pi/\k$-periodic wave train of \eqref{E:main} and \eqref{def:c1} is modulationally unstable, provided that $\Delta(\k)<0$. It is spectrally stable to square integrable perturbations in the vicinity of the origin in $\mathbb{C}$ otherwise. A change in sign of $\Delta(\k)$ and thus, in stability occurs when one of the factors $i_j$'s, $j=1,2,3,4$ vanishes. Notice that for a fixed $T$, all these factors explicitly depend on the wave number $\k$, the phase velocity $\c(\k)$, and the group velocity $(\k\c(\k))'$. Therefore, the vanishing of each of the factor is associated with some resonance in the wave (see \cite{HP2}). Specifically,
\begin{itemize}
\item[(R1)] $i_1(\k)$ is derivative of the group velocity and therefore, if $i_1(\k)=0$ at some $\k$, the group velocity achieves an extremum at the wave number $\k$; 
\item[(R2)] $i_2(\k)$ is the difference between the group velocity and the phase velocity in the long wave limit as $\k\to0$, that is, $\pm \c(0)=\pm 1$ and therefore, if $i_2(\k)=0$  at some $\k$; it results in the ``resonance of short and long waves;"
\item[(R3)] $i_3(\k)$ is the difference between the phase velocities of the fundamental mode, $\pm \c(\k)$ and second harmonic, $\pm \c(2\k)$ and therefore, $i_3(\k)=0$ at some $\k$ implies ``second harmonic resonance;"
\item[(R4)] $i_4(\k)$ is the only factor which captures the nonlinearity, and we expect $i_4(\k)$ to vanish when dispersion effects balance the nonlinear effects. 
\end{itemize}

For $T=0$, $i_1(\k)<0$ and $i_2(\k)<0$ for any $\k>0$ while $i_3(\k)>0$ for any $\k>0$. A numerical evaluation of \eqref{def:i4} reveals a unique root $\k_c=1.008\dots$ of $i_4$ over the interval $(0,\infty)$ such that $i_4(\k)>0$ if $0<\k<\k_c$ and it is negative if $\k_c<\k<\infty$. Therefore, for $T=0$, a sufficiently small $2\pi/\k$-periodic traveling wave of \eqref{E:main} and \eqref{def:c} is modulationally unstable if $\k>\k_c$. It is modulationally stable if $0<\k<\k_c$. 

For $T>0$, we describe the modulational instability through the diagram Figure~\ref{fig1:d}. In $\k$-$\k\sqrt{T}$ plane, four curves are corresponding to each mechanism splitting the plane into three regions of stability and three regions of instability. Any fixed $T>0$ corresponds to a line passing through the origin of slope $\sqrt{T}$. For $0<T<1/3$, the line crosses all the curves producing three intervals of stable wave numbers and three intervals of unstable wave numbers. Therefore, for $0<T<1/3$, all the four mechanisms (R1) to (R4) contribute towards modulational instability.

On the other hand, for $T>1/3$, the line through the origin only crosses the Curve~4 corresponding to $i_4(\k)=0$, see Figure~\ref{fig1:d}. In this case, the modulational instability is caused only by the mechanism (R4) similar to the case $T=0$. For every $T>1/3$, there is a critical wave number $\k_c(T)$ such that a sufficiently small $2\pi/\k$-periodic traveling wave of \eqref{E:main} and \eqref{def:c1} is modulationally unstable if $\k>\k_c(T)$. Moreover, $\lim_{T\to \infty} \k_c(T) = \infty$. The result becomes inconclusive for $T=1/3$.    

%%%%%%%%%%%%%%%%%%%%%%%%%%%%%%%%%%%%%%%%%%%%%%%%%%
%%%%%%%%%%%%%%%%%%%%%%%%%%%%%%%%%%%%%%%%%%%%%%%%%%
%%%%%%%%%%%%%%%%%%%%%%%%%%%%%%%%%%%%%%%%%%%%%%%%%%
\section{Comparison with other models}\label{sec:T>0}
%%%%%%%%%%%%%%%%%%%%%%%%%%%%%%%%%%%%%%%%%%%%%%%%%%
%%%%%%%%%%%%%%%%%%%%%%%%%%%%%%%%%%%%%%%%%%%%%%%%%%
%%%%%%%%%%%%%%%%%%%%%%%%%%%%%%%%%%%%%%%%%%%%%%%%%%

%%%%%%%%%%%%%%%%%%%%%%%%%%%%%%%%%%%%%%%%%%%%%%%%%%
%%%%%%%%%%%%%%%%%%%%%%%%%%%%%%%%%%%%%%%%%%%%%%%%%%
%%%%%%%%%%%%%%%%%%%%%%%%%%%%%%%%%%%%%%%%%%%%%%%%%%
\begin{table}
\centering
\begin{tabular}{r|c}
\textbf{Model} & $\boldsymbol{i_4(\k)}$ \\ \hline
\textbf{Whitham} & $(2 i_3^-+i_2^-)(\k)$ \\ \hline
\textbf{FDCH} & $\begin{aligned} \Big( 3i_2^- - i_2^- i_3^- + 6i_3^-  -\tfrac{1}{12} \k^2(57i_2^- +34 i_3^-) +\tfrac{1}{108}\k^4(198i_2^- + 35i_3^-)\Big)(\k) \end{aligned}$\\ \hline
\textbf{FDSW-I} & $\begin{aligned}3\c^2(\k)+5\c^4(\k)&-  2\c^2(2\k)  (\c^2(\k)+2)+18\k\c^3(\k)\c'(\k)\\&+\k^2(\c')^2(\k)(5\c^2(\k)+4\c^2(2\k))\end{aligned}$ \\ \hline
\end{tabular}
\caption{$i_4(\k)$ for different full-dispersion shallow water models.}
\label{table1}
\end{table}
%%%%%%%%%%%%%%%%%%%%%%%%%%%%%%%%%%%%%%%%%%%%%%%%%%
%%%%%%%%%%%%%%%%%%%%%%%%%%%%%%%%%%%%%%%%%%%%%%%%%%
%%%%%%%%%%%%%%%%%%%%%%%%%%%%%%%%%%%%%%%%%%%%%%%%%%

The index formula \eqref{def:ind} for FDSW-II is similar to the index formula for the FDSW-I model (see \cite{HP2}) except $i_4(\k)$, see Table~\ref{table1}. It is a manifestation of the fact that $i_1, i_2$, and $i_3$ are completely given in terms of the dispersion and independent of the nonlinearity of the equation. Since, both the FDSW-I and FDSW-II models have same dispersion, $i_1, i_2$, and $i_3$ are same for both the models. The difference in $i_4$ reveals that how different these two models are as far as the modulational instability is concerned. The index formulas for the Whitham and FDCH equations (see \cite{HJ2} and \cite{HP3}) are different from \eqref{def:ind} not only in $i_4$, see Table~\ref{table1}, but also in $i_2$ and $i_3$. The factors $i_2$ and $i_3$ are replaced by $i_2^-$ and $i_3^-$ in index formulas for the Whitham and FDCH equations, where
\[
i_2(\k) = ((\k \c(\k))'-1) ((\k \c(\k))'+1) =: i_2^-(\k) i_2^+(\k) \quad \text{and} 
\]  
\[
\quad i_3(\k)=(\c(\k)-\c(2\k)) (\c(\k)+\c(2\k)) =: i_3^-(\k) i_3^+(\k). 
\]
Due to their unidirectionality, index formulas for the Whitham and FDCH equations only capture the resonances of the positive branch of the dispersion of full water wave equations and miss the resonances between positive and negative branches of dispersion. Such resonances are present in the FDSW-I, and FDSW-II models as their dispersion have both branches of the dispersion of water waves.

All the models, Whitham, FDCH, FDSW-I and FDSW-II qualitatively exhibit Benjamin-Feir instability of a Stokes wave with same mechanism (R4); see \cite{HJ2,HP2,HP3}.  Although all these models have the dispersion of the full water wave problem, they are approximate models, and it should not come as a surprise that each of them predicts the Benjamin-Feir instability with a different critical wave number (see Table~\ref{table}), not equal to $1.363\dots$.

%%%%%%%%%%%%%%%%%%%%%%%%%%%%%%%%%%%%%%%%%%%%%%%%%%
%%%%%%%%%%%%%%%%%%%%%%%%%%%%%%%%%%%%%%%%%%%%%%%%%%
%%%%%%%%%%%%%%%%%%%%%%%%%%%%%%%%%%%%%%%%%%%%%%%%%%
\begin{table}
\centering
\begin{tabular}{r|c}
\textbf{Model} & \textbf{Critical wave number} \\ \hline
\textbf{Full water wave} & $1.363$ \\ \hline
\textbf{Whitham} & $1.363-0.217$ \\ \hline
\textbf{FDCH} & $1.363+0.057$ \\ \hline
\textbf{FDSW-I} & $1.363+0.247$ \\ \hline
\textbf{FDSW-II} & $1.363-0.355$ \\
\end{tabular}
\caption{The value of critical wave number for different full-dispersion shallow water models.}
\label{table}
\end{table}
%%%%%%%%%%%%%%%%%%%%%%%%%%%%%%%%%%%%%%%%%%%%%%%%%%
%%%%%%%%%%%%%%%%%%%%%%%%%%%%%%%%%%%%%%%%%%%%%%%%%%
%%%%%%%%%%%%%%%%%%%%%%%%%%%%%%%%%%%%%%%%%%%%%%%%%%

\begin{figure}[ht] 
  \begin{subfigure}[b]{0.49\linewidth}
    \centering
    \includegraphics[scale=0.49]{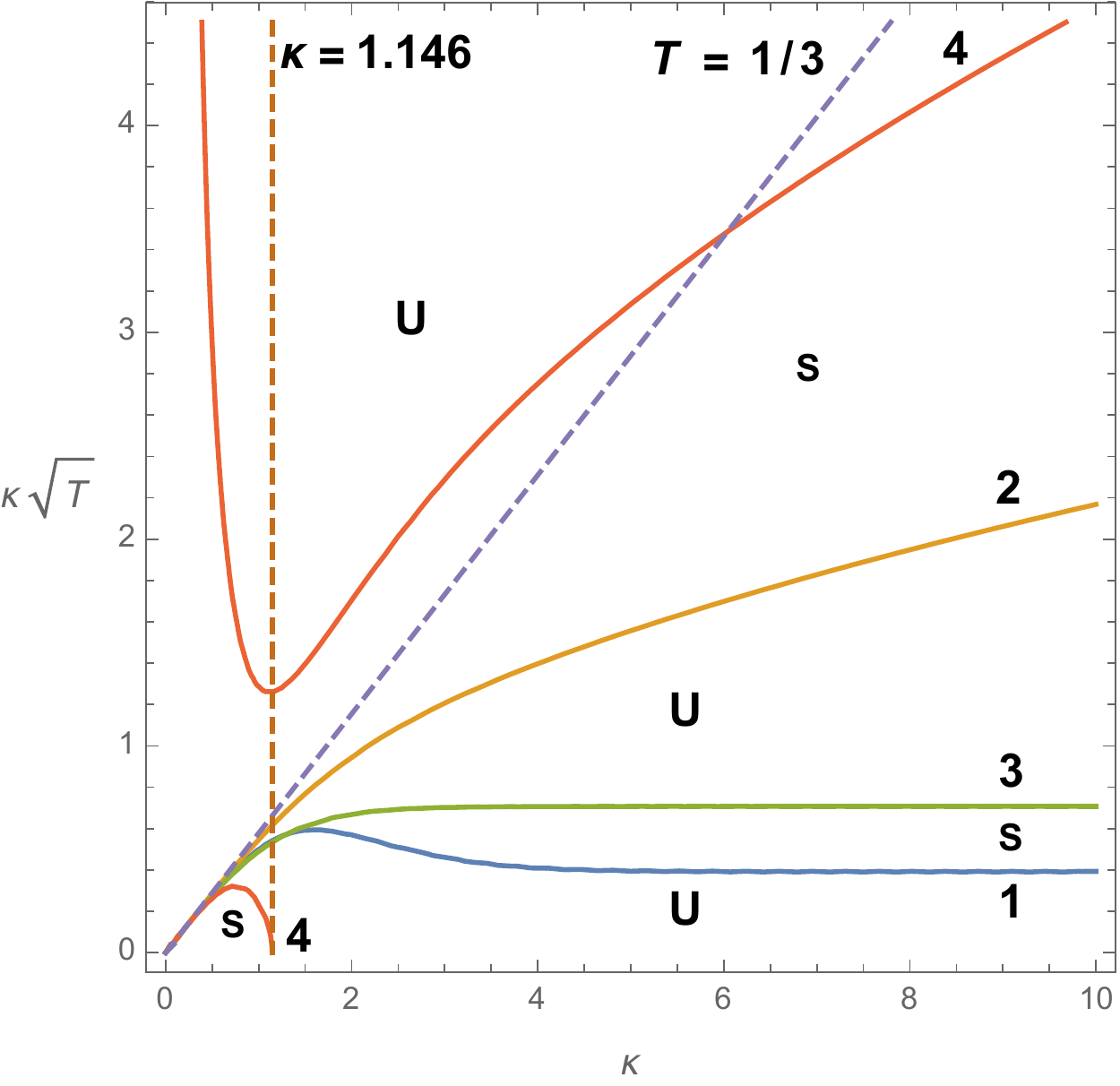} 
    \caption{Whitham} 
    \label{fig1:a} 
    \vspace{1ex}
  \end{subfigure} 
   \begin{subfigure}[b]{0.5\linewidth}
    \centering
    \includegraphics[scale=0.39]{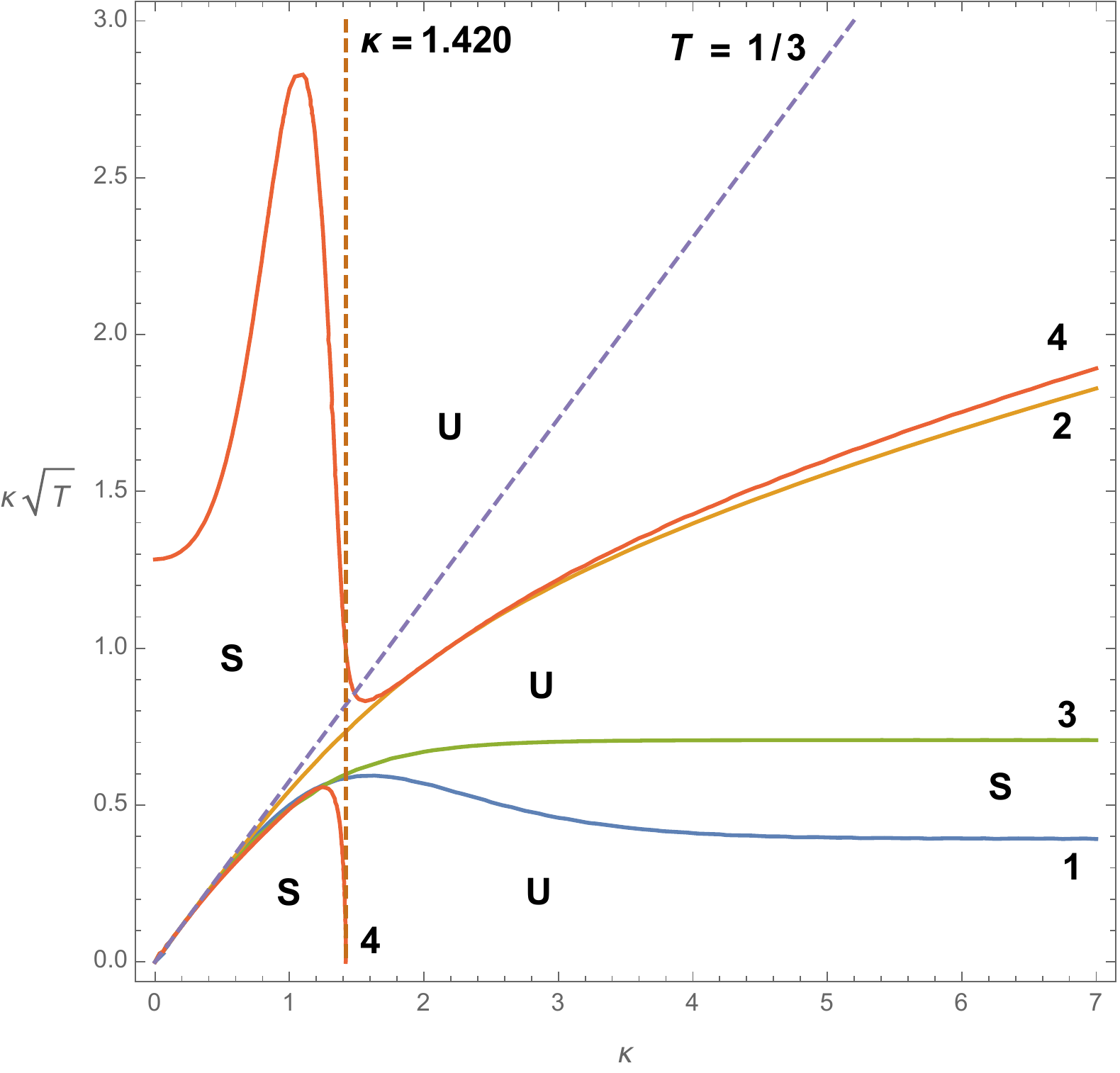} 
    \caption{FDCH} 
    \label{fig1:b} 
    \vspace{1ex}
  \end{subfigure} %%
  \begin{subfigure}[b]{0.5\linewidth}
    \centering
    \includegraphics[scale=0.49]{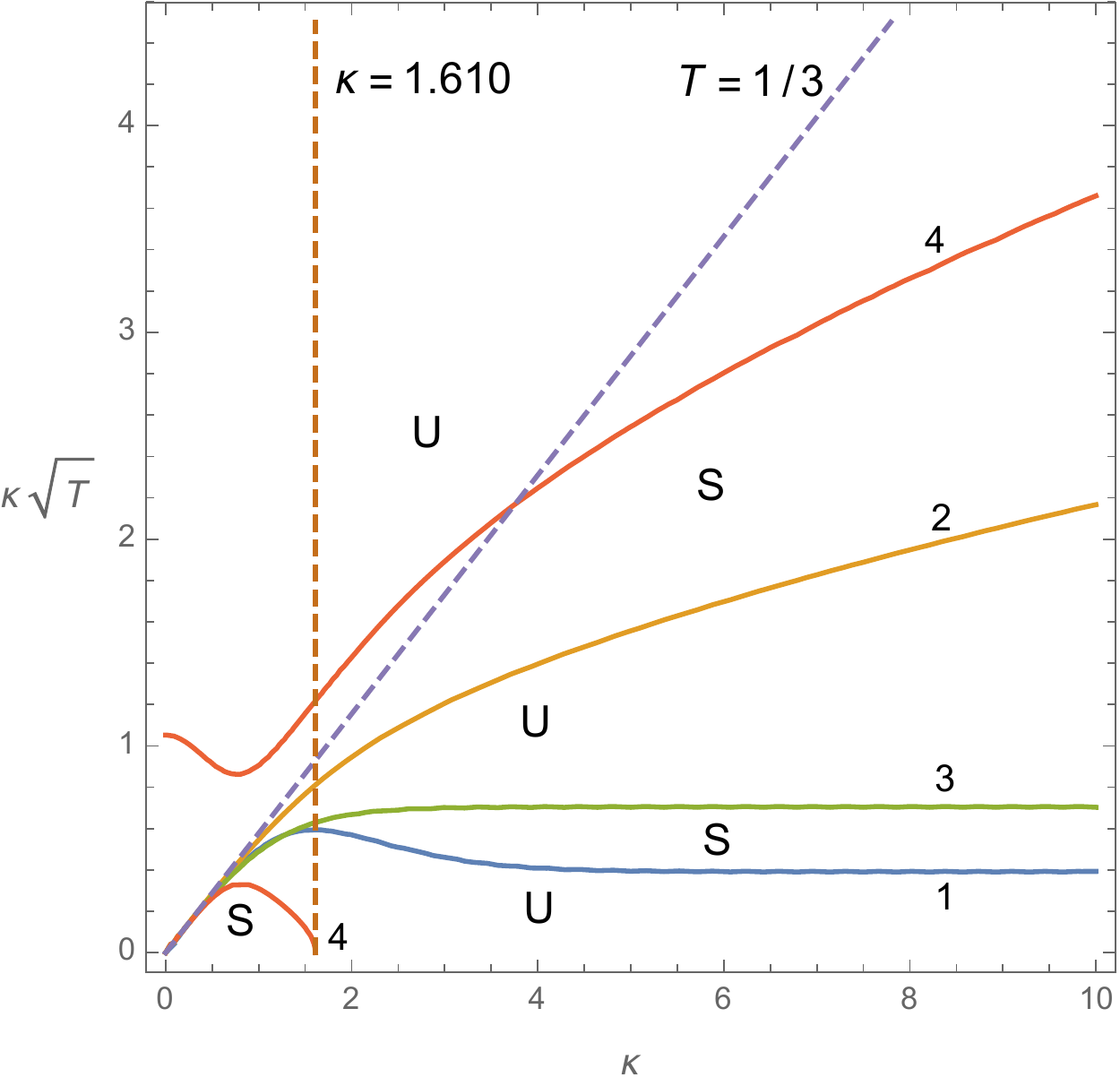} 
    \caption{FDSW-I} 
    \label{fig1:c} 
  \end{subfigure}%%
  \begin{subfigure}[b]{0.5\linewidth}
    \centering
    \includegraphics[scale=0.49]{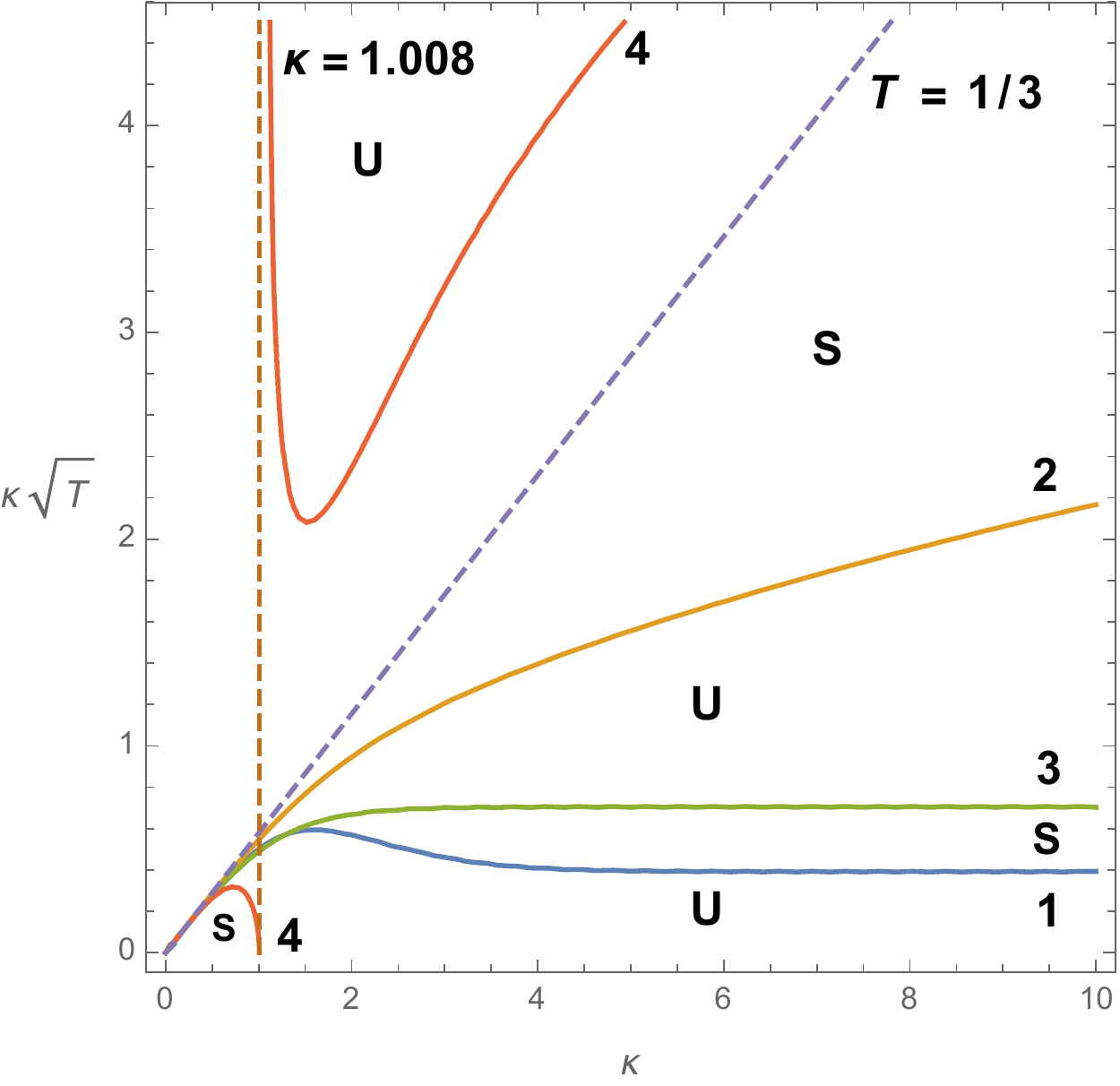} 
    \caption{FDSW-II} 
    \label{fig1:d} 
  \end{subfigure} 
   \begin{subfigure}[b]{0.7\linewidth}
    \centering
    \includegraphics[width=8cm, height=6cm]{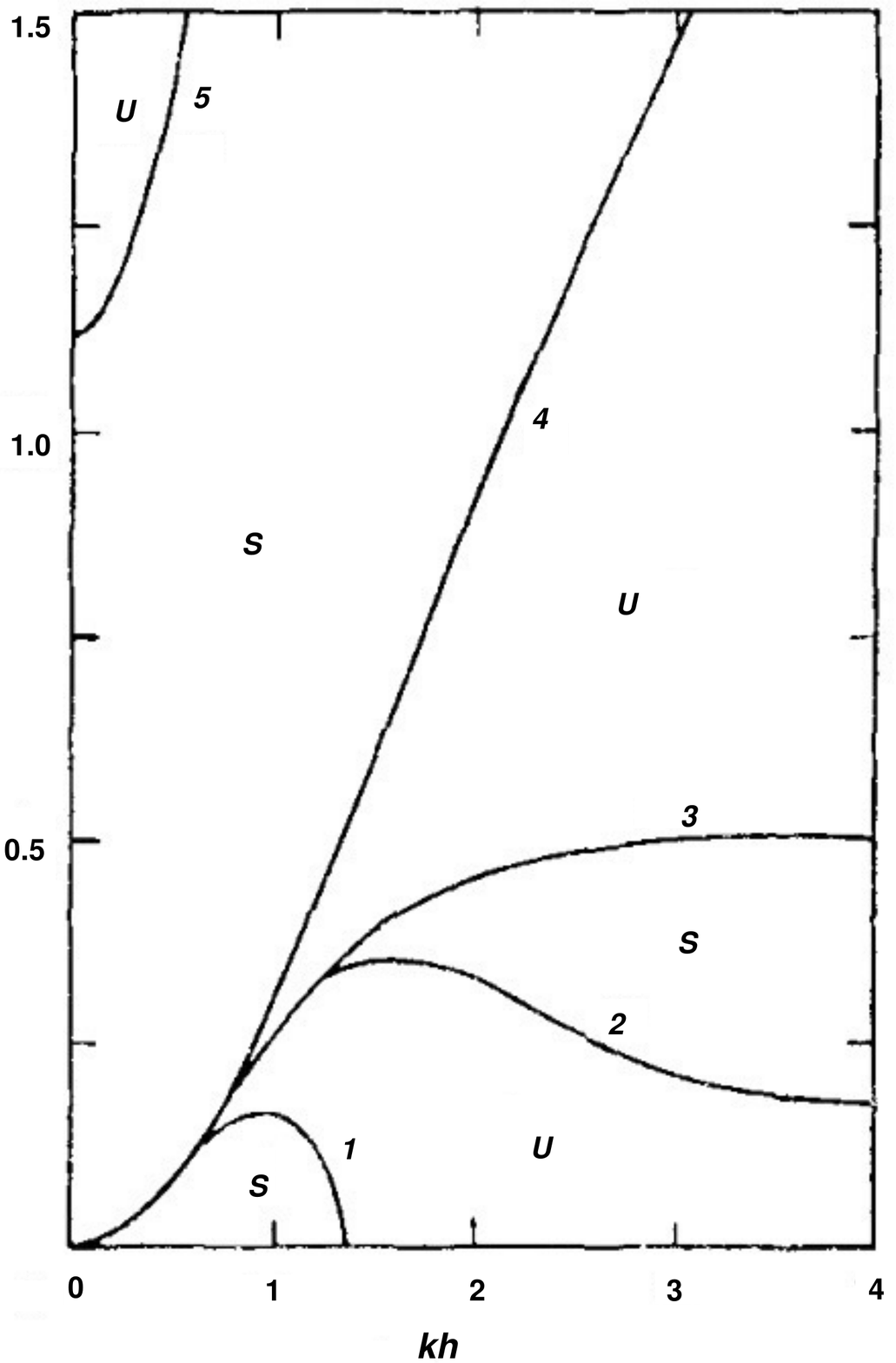} 
    \caption{Full water wave} 
    \label{fig1:e} 
  \end{subfigure}%% 
   \caption{Stability diagram for sufficiently small, periodic wave trains of models indicated. ``S" and ``U" denote stable and unstable regions. In Figures~\ref{fig1:a}-\ref{fig1:d}, solid curves represent roots of the modulational instability index and are labeled according to their mechanism. Figure~\ref{fig1:a}, \ref{fig1:b} and \ref{fig1:c} are adapted from \cite{HJ3}, \cite{HP3} and \cite{HP2} respectively. Figure~\ref{fig1:e} is taken from \cite{DR}.}
  \label{fig1} 
\end{figure}

The effects of surface tension on modulational instability in all these models along with the full water wave problem have been compared in Figure~\ref{fig1}. The diagrams corresponding to model equations, Figure~\ref{fig1:a},\ref{fig1:b},\ref{fig1:c},\ref{fig1:d}, contain four curves corresponding to each mechanism from (R1) to (R4). The diagram corresponding to the physical problem, Figure~\ref{fig1:e}, has five curves and by a direct comparison with the model equations, it can be deduced that Curves~2, 3 and 4 are coming from mechanisms (R1), (R2) and (R3) respectively since the full water wave problem shares dispersion with all these models. Moreover, Curves~1 and 5 of Figure~\ref{fig1:e} can be results of the interaction between dispersion and nonlinearity of the full water wave problem, like other models.  

In Figure~\ref{fig1}, a fixed $T > 0$ corresponds to a line passing through the origin. For small surface tensions, in the physical problem, the wave numbers are divided into three intervals of stability and three intervals of instability, see Figure~\ref{fig1:e}. All the models agree with the physical problem for small surface tensions, more precisely, for $0<T<1/3$. 

The physical problem reveals that for sufficiently large surface tension, the stability changes to instability only once about a critical wave number much like Benjamin-Feir instability for $T=0$, see Figure~\ref{fig1:e}. In all the models, for $T>1/3$, there is a critical wave number $\k_c(T)$ about which the stability changes to instability and therefore, all the models agree qualitatively with the physical problem. The difference in models arises when we look at $\lim_{T \to \infty} \k_c(T)$. The physical problem suggests that this limit is finite and approximately equal to $1.121$.  As we discussed in Section~\ref{sec:MI}, $\lim_{T \to \infty} \k_c(T)$ diverges for FDSW-II model. In other words, all sufficiently small periodic traveling waves of \eqref{E:main} and \eqref{def:c1} are modulationally stable in the large surface tension limit, which is unphysical as suggested by the physical problem. On the other hand, for FDSW-I model, $\lim_{T \to \infty} \k_c(T)\approx 1.054$. Therefore, in the large surface tension limit, the FDSW-I model explains the effects of surface tension similar to the physical problem. For the Whitham equation,  $\lim_{T \to \infty} \k_c(T)$ diverges, and it fails to explain the effects in large surface tension limit. For the FDCH equation, $\lim_{T \to \infty} \k_c(T)\approx 1.283$ and it offers an improvement over the Whitham equation.      

The comparative study suggests that although both FDSW-I and FDSW-II are bi-directional shallow water models extending nonlinear shallow water equations to include full-dispersion of water waves, the FDSW-I is a better model as far as modulational instability is concerned.      

%%%%%%%%%%%%%%%%%%%%%%%%%%%%%%%%%%%%%%%%%%%%%%%%%%
%%%%%%%%%%%%%%%%%%%%%%%%%%%%%%%%%%%%%%%%%%%%%%%%%%
%%%%%%%%%%%%%%%%%%%%%%%%%%%%%%%%%%%%%%%%%%%%%%%%%%
%%%%%%%%%%%%%%%%%%%%%%%%%%%%%%%%%%%%%%%%%%%%%%%%%%
%%%%%%%%%%%%%%%%%%%%%%%%%%%%%%%%%%%%%%%%%%%%%%%%%%
%%%%%%%%%%%%%%%%%%%%%%%%%%%%%%%%%%%%%%%%%%%%%%%%%%

%%%%%%%%%%%%%%%%%%%%%%%%%%%%%%%%%%%%%%%%%%%%%%%%%%
%%%%%%%%%%%%%%%%%%%%%%%%%%%%%%%%%%%%%%%%%%%%%%%%%%
%%%%%%%%%%%%%%%%%%%%%%%%%%%%%%%%%%%%%%%%%%%%%%%%%%
\subsection*{Acknowledgements}

The author thanks Vera Mikyoung Hur for helpful discussions. 

\bibliographystyle{amsalpha}
\bibliography{stability.bib}

\end{document}